\documentclass[]{interact}
\usepackage{amsmath}
\usepackage{amssymb}
\usepackage{epstopdf}
\usepackage[caption=false]{subfig}
\usepackage[numbers,square,sort&compress]{natbib}
\usepackage{lmodern}
\usepackage{url}
\usepackage[hidelinks]{hyperref}

\newcommand{\E}{\mathbb{E}}

\theoremstyle{plain}
\newtheorem{theorem}{Theorem}[section]
\newtheorem{lemma}[theorem]{Lemma}
\newtheorem{corollary}[theorem]{Corollary}

\theoremstyle{definition}
\newtheorem{definition}[theorem]{Definition}

\theoremstyle{remark}

\begin{document}

\title{Quantum Adiabatic Theorem under Wiener Noise: Maximal Mean-Square Estimates}

\author{\name{A.M.H.H. Caldera\textsuperscript{1 $\ast$}\thanks{
\textsuperscript{$\ast$}rok.harshahiranya@gmail.com},
R.P.K. De Silva\textsuperscript{2 $\dagger$}\thanks{
\textsuperscript{$\dagger$}kaushika.desilva@uwa.edu.au},
K.K.W.A.S. Kumara\textsuperscript{1 $\ddagger$}\thanks{
\textsuperscript{$\ddagger$}sarath@sjp.ac.lk}
and B.P.W. Fernando\textsuperscript{1 \S}\thanks{
\textsuperscript{\S}Corresponding author: bpw@sci.sjp.ac.lk}}
\affil{\textsuperscript{1}Department of Mathematics,
University of Sri Jayewardenepura, Sri Lanka;
\textsuperscript{2}Centre for Quantum Information,
Simulation and Algorithms, University of Western Australia, Australia}}

\maketitle

\begin{abstract}
We study finite-dimensional stochastic Schrödinger equations with Hamiltonian-proportional Wiener noise. For deterministic \(C^2\)-Hamiltonians whose eigenvalues have constant multiplicities and do not cross, we construct a unitary stochastic adiabatic intertwiner driven by the same Wiener process as the exact stochastic propagator. A commutator identity, Itô calculus, and Doob’s maximal inequality bound their operator-norm difference on a selected initial eigenspace. Under uniform coefficient bounds, physical noise scaling \(\alpha=\sigma T^{-\gamma}\), with fixed \(\sigma\geq0\) and \(\gamma\geq0\), yields maximal mean-square error \(O(T^{-2}+\sigma^4T^{-4\gamma}+\sigma^2T^{-(2\gamma+1)})\), establishing convergence for every \(\gamma>0\). For fixed nonzero physical noise amplitude, the estimate remains inconclusive.
\end{abstract}


\section{Introduction}
The quantum adiabatic theorem has its origins in the 1928 analysis of slowly varying Hamiltonians by Born and Fock~\cite{BornFock1928}. In 1950, Kato~\cite{Kato1950} developed an operator-theoretic formulation based on spectral projections, accommodating isolated
eigenvalues of finite multiplicity. For a fixed, sufficiently regular path of self-adjoint Hamiltonians \(H_T(s)\), \(s\in[0,1]\), increasing the total evolution time \(T\) in the driving \(H(t/T)\), \(0\leq t\leq T\), slows the physical variation. When the relevant eigenvalue remains uniformly separated from the rest of the spectrum, the standard adiabatic theorem ensures that a state initially in its eigenspace remains close to the corresponding instantaneous eigenspace, with the distance in Hilbert space norm tending to zero uniformly in \(s\) as \(T\to\infty\). The associated adiabatic approximation includes the dynamical phase and geometric evolution within that eigenspace. Avron and Elgart~\cite{AvronElgart1999} later established adiabatic convergence in certain gapless settings for finite-rank spectral projections.

Our main result, Theorem~\ref{thm:stochastic version quantum adiabatic}, adapts the deterministic argument in Section~8.2 and Appendix~G of ~\cite{Scherer2019}. After rescaling physical time to \(s\in[0,1]\), the exact propagator \(U_T(s)\) satisfies the Schrödinger equation
\[
i\dot U_T(s)=T H_T(s)U_T(s).
\]
For a differentiable spectral projection \(P_j(s)\) satisfying
\([H_T(s),P_j(s)]=0\), the corresponding adiabatic intertwiner \(U_{A,j}(s)\) is defined by
\[
i\dot U_{A,j}(s)
=
\left(TH_T(s)+i[\dot P_j(s),P_j(s)]\right)U_{A,j}(s).
\]
It obeys the exact intertwining relation
\[
P_j(s)U_{A,j}(s)=U_{A,j}(s)P_j(0).
\]
For a uniformly isolated eigenvalue branch of constant multiplicity, we can define a resolvent-based operator which gives us the commutator identity that we use to integrate by parts when comparing \(U_T\) and \(U_{A,j}\). Under uniform spectral gap and \(C^2\)-regularity bounds, this argument gives us
\[
\sup_{0\le s\le 1}\left\|(U_T(s)-U_{A,j}(s))P_j(0)\right\|
=
O(T^{-1}).
\]
The theorem~G.15 of ~\cite{Scherer2019} builds on the quantitative operator estimates of Jansen, Ruskai and Seiler~\cite{JansenRuskaiSeiler2007}.

Random fluctuations in the control fields or environmental parameters introduce contributions that this deterministic estimate cannot directly control. Wiener noise provides an idealized model of continuous random forcing, including stochastic phase fluctuations. For a fixed Hamiltonian path, increasing the total evolution time \(T\) slows the evolution and allows the noise to act for longer. Therefore, the question is under what conditions and for what scaling of the physical noise amplitude with the total evolution time \(T\) the exact stochastic evolution approaches the appropriate stochastic adiabatic evolution.

The existing approaches address several aspects of this problem. Barchielli, Pellegrini and Petruccione~\cite{Barchielli2010} model non-Markovian effects in the stochastic Schrödinger equation using colored noise, while Gaitan~\cite{Gaitan2006} studies adiabatic search performance under random Zeeman-coupled perturbations. Chinthila et al.~\cite{Chinthila2025} simulate an adiabatic algorithm for Exact Cover-3 using a semimartingale approximation to fractional Wiener process. They estimate the probability that a measurement at the final time yields a valid solution to the Exact Cover-3 instance. Avron, Fraas, Graf and Grech~\cite{AvronFraasGrafGrech2012} develop an adiabatic theory for generators of contracting evolutions, including Lindblad generators. Fraas~\cite[Theorem~5]{Fraas2014} derives a stochastic adiabatic expansion and determines the leading-order tunneling statistics for a dephasing model with simple eigenvalues. In this model, the stochastic trajectories acquire random dynamical phases, and the ensemble-averaged density matrix satisfies a Lindblad equation.

Our study gives an explicit maximal mean-square estimate for the difference between the exact stochastic propagator \(U_T(s)\) and stochastic adiabatic intertwiner \(U_{A,j}(s)\) on the selected initial eigenspace \(P_j(0)\mathcal H\), both being driven by the same real-valued Wiener process. We assume that the finite-dimensional Hamiltonian is twice continuously differentiable on \([0,1]\), that every eigenvalue has constant multiplicity, and that distinct eigenvalue branches do not cross. We consider a selected eigenvalue separated from the rest of the spectrum by a uniform positive gap on this interval. The diffusion operator \(B(t)\) is proportional to the Hamiltonian \(H(t)\), so it commutes with the instantaneous orthogonal spectral projection. This results in a unitary stochastic adiabatic intertwiner that preserves the deterministic intertwining structure.

For the physical noise scaling \(\alpha=\sigma T^{-\gamma}\), with fixed \(\sigma\), the estimate establishes stochastic adiabatic convergence for every \(\gamma>0\), provided its spectral and derivative coefficients remain bounded independently of \(T\). At \(\gamma=\tfrac{1}{2}\), the maximal mean-square error \(\mathcal E_T\) is \(O(T^{-2})\), corresponding to an \(O(T^{-1})\) root-mean-square maximal error \(\mathcal E_T^{1/2}\) and retaining the usual deterministic adiabatic order. For fixed nonzero physical noise amplitude, \(\gamma=0\), the upper bound contains a term that does not vanish as the total evolution time \(T\) increases. Therefore, it leaves convergence unresolved and does not prove the breakdown of adiabaticity. 

\section{Preliminaries}
Throughout this paper, \(\|\cdot\|\) denotes the Hilbert space norm for
vectors and the operator norm for operators. For the Hilbert space \(\mathcal H\) of fixed finite dimension and \(A\in\mathcal B(\mathcal H)\),
\[
\|A\|=\sup_{\|x\|=1}\|Ax\|,
\qquad
x\in\mathcal H.
\]
Here $\mathcal{B}(\mathcal{H})$ denotes the set of bounded linear operators on the Hilbert space $\mathcal{H}$.
\begin{definition}[Time-rescaled operators]
\label{def:time rescaled operators}
Let $t_{\mathrm{ini}}$ and $t_{\mathrm{fin}} (> t_{\mathrm{ini}})$ denote the initial and final times of the evolution. We define the total elapsed time by
\begin{equation}
T = t_{\mathrm{fin}} - t_{\mathrm{ini}}.
\label{eq:1}
\end{equation}
To reformulate the dynamics on a fixed parameter interval, we introduce the dimensionless rescaled time variable $s \in [0,1]$. The corresponding rescaled Hamiltonian $H_T$ is defined by
\begin{align}
H_T : [0,1]&\to\mathcal{B}_{\mathrm{sa}}(\mathcal{H}),\notag\\
\quad
H_T(s)&=H\bigl(t_{\mathrm{ini}} + sT\bigr),
\qquad s\in[0,1],
\label{eq:2}
\end{align}
where $\mathcal{B}_{\mathrm{sa}}(\mathcal{H})$ denotes the space of bounded self-adjoint operators on the Hilbert space $\mathcal{H}$.

\medskip
\noindent
Similarly, we define the time-rescaled evolution operator $U_T$ by
\begin{align}
U_T : [0,1]&\to\mathcal{B}(\mathcal{H}),\notag\\
\quad
U_T(s)&=U\bigl(t_{\mathrm{ini}} + sT,\, t_{\mathrm{ini}}\bigr),
\qquad s\in[0,1],
\label{eq:3}
\end{align}
where $U(t,t_{\mathrm{ini}})$ denotes the evolution operator associated with the time-dependent Hamiltonian $H(t)$.
This time rescaling transfers the original dynamics from the physical time interval $[t_{\mathrm{ini}},t_{\mathrm{fin}}]$ to the fixed interval $[0,1]$.
\end{definition}

\medskip
\noindent
As an immediate consequence of \eqref{eq:2}, for every $t \in [t_{\mathrm{ini}},t_{\mathrm{fin}}]$ one has
\begin{equation}
H_T\!\left(\frac{t-t_{\mathrm{ini}}}{T}\right) = H(t).
\label{eq:4}
\end{equation}
In particular, evaluating the rescaled Hamiltonian at the endpoints yields
\begin{align*}
H_T(0)
&= H(t_{\mathrm{ini}})\notag\\
&= H_{\mathrm{ini}}\notag\\
H_T(1)
&= H(t_{\mathrm{fin}})\notag\\
&= H_{\mathrm{fin}}.
\end{align*}
This formulation is especially useful in the stochastic adiabatic setting, where one seeks to distinguish rescaled time \(s\) from additional random effects in the dynamics.
\begin{lemma}[Time-rescaled evolution]
\label{lem:time rescaled evolution}
The time-rescaled evolution operator $U_T(\cdot)$ satisfies
\begin{equation}
i\frac{d}{ds}U_T(s)=T\,H_T(s)\,U_T(s),
\qquad s\in[0,1],
\label{eq:5}
\end{equation}
together with the boundary conditions
\begin{equation}
U_T(0)=I,
\qquad
U_T(1)=U(t_{\mathrm{fin}},t_{\mathrm{ini}}).
\label{eq:6}
\end{equation}
\end{lemma}
\begin{proof}
See~\cite[Lemma~8.2, p.~407]{Scherer2019}.
\end{proof}
\subsection{Adiabatic Assumption}
Let
\begin{align*}
H_T : [0,1] &\to \mathcal{B}_{\mathrm{sa}}(\mathcal{H})\notag\\
\quad
s &\mapsto H_T(s)
\end{align*}
be the time-rescaled Hamiltonian on the Hilbert space $\mathcal{H}$. We assume that the following conditions are satisfied by $H_T(s)$:

\begin{enumerate}
    \item[(i)] \textit{Regularity:}
    The map $s \mapsto H_T(s)$ is at least twice continuously differentiable on the closed interval $[0,1]$, that is,
    \[
    H_T \in C^2\bigl([0,1], \mathcal{B}_{\mathrm{sa}}(\mathcal{H})\bigr).
    \]

    \item[(ii)] \textit{Purely discrete spectrum with constant multiplicities:}
    For every $s \in [0,1]$, the operator $H_T(s)$ has purely discrete spectrum
    \[
    \sigma\bigl(H_T(s)\bigr)=\{E_j(s)\mid j\in \mathbf{I}\subset \mathbb{N}_0\}.
    \]
    Each eigenvalue $E_j(s)$ may be degenerate, with multiplicity $d_j(s)$, and these multiplicities are assumed to remain constant on $[0,1]$, i.e.
    \[
    d_j(s)=d_j
    \qquad \text{for all} \ s\in [0,1], \quad j\in \mathbf{I}.
    \]

    \item[(iii)] \textit{Absence of level crossings:}
    The eigenvalues of $H_T(s)$ do not cross on $[0,1]$. More precisely, they can be ordered so that for every $s\in [0,1]$,
    \[
    E_0(s)<E_1(s)<\cdots<E_j(s)<E_{j+1}(s)<\cdots.
    \]
\end{enumerate}
The assumption (i) provides the smoothness required in the derivation of the adiabatic approximation and in the control of the associated error terms. The assumption (ii) ensures that the relevant spectral subspaces are well defined throughout the evolution. The assumption (iii) excludes eigenvalue crossings and guarantees that distinct spectral branches remain separated, which is a basic requirement in the standard formulation of the quantum adiabatic theorem.
\begin{definition}[The resolvent]
\label{def:resolvent}
The resolvent $R_z(A)$ of an self-adjoint bounded operator \(A\) on a Hilbert space $\mathcal{H}$ is
defined as
\begin{align}
R_{(.)}(A) : \mathbb{C}\setminus\sigma(A) &\to \mathcal{L}(\mathcal{H})\notag\\
\quad
z &\mapsto R_z(A)=(A-zI)^{-1}.
\label{eq:7}
\end{align}
Here, $\mathcal{L}(\mathcal{H})$ denotes the space of linear operators on the Hilbert space $\mathcal{H}$.
\end{definition}
\begin{definition}[The resolvent-based operator]
\label{def:the resolvent-based operator}
Let $H_T$ be an operator-valued function satisfying the Adiabatic Assumption. Moreover, let $\gamma_j(s)$ be any closed counter-clockwise curve in $\mathbb{C}\setminus\sigma(H_T(s))$ that encloses solely the eigenvalue \(E_j(s)\) from \(\sigma(H_T(s))\). For any self-adjoint bounded operator \(A\in C^1([0,1])\),
\begin{align*}
A : [0,1] &\to \mathcal{B}(\mathcal{H})\notag\\
\quad
s &\mapsto A(s)
\end{align*}
and \(j\in \mathbf{I}\) define
\begin{align*}
X_j[A] : [0,1] &\to \mathcal{L}(\mathcal{H})\notag\\
\quad
s &\mapsto X_j[A](s)
\end{align*}
by
\begin{equation}
X_j[A](s)
=
\frac{1}{2{\pi}i}\int_{\gamma_j(s)}
R_z(H_T(s))A(s)R_z(H_T(s))\,dz,
\qquad s\in[0,1].
\label{eq:8}
\end{equation}
\end{definition}
\begin{definition}[Gap function]
\label{def:gap function}
Let \(H_T\) be an operator-valued function satisfying the Adiabatic Assumption. For each eigenvalue \(E_j(s)\) of \(H_T(s)\), the associated gap function \(g_j:[0,1]\to \mathbb{R}\) is defined by
\begin{equation}
g_j(s)=\min_{k\in \mathbf{I}\setminus\{j\}} \lvert E_k(s)-E_j(s)\rvert,
\qquad s\in[0,1].
\label{eq:9}
\end{equation}
\end{definition}
\begin{lemma}[Spectral projections]
\label{lem:spectral-projections}
Let \(H_T\) be an operator-valued function satisfying the Adiabatic Assumption and let \(\{P_j(s)\}_{j\in \mathbf{I}}\) denote the family of orthogonal spectral projections onto the eigenspaces of \(H_T(s)\). For each \(j\in \mathbf{I}\) and \(s\in[0,1]\), one has
\begin{equation}
\dot{P}_j(s)=X_j\!\left[\dot{H}_T\right](s)
\label{eq:10}
\end{equation}
and
\begin{equation}
\left[H_T(s),\,X_j[\dot{P}_j](s)\right]
=
-\left[\dot{P}_j(s),\,P_j(s)\right],
\qquad s\in[0,1].
\label{eq:11}
\end{equation}
\end{lemma}
\begin{proof}
See~\cite[Lemma~G.5, p.~624]{Scherer2019}.
\end{proof}
\begin{definition}[Adiabatic generator]
\label{def:adiabatic-generator}
Let \(H_T\) be an operator-valued function satisfying the Adiabatic Assumption. 
For each \(j\in \mathbf{I}\), the associated adiabatic generator \(H_{A,j}(s)\) is defined by
\begin{align}
H_{A,j}:[0,1]&\to\mathcal{B}_{\mathrm{sa}}(\mathcal{H}),\notag\\
\quad
H_{A,j}(s)&=T\,H_T(s)+i\left[\dot{P}_j(s),\,P_j(s)\right],
\qquad s\in[0,1].
\label{eq:12}
\end{align}
\end{definition}
\begin{definition}[Adiabatic intertwiner]
\label{def:adiabatic-intertwiner}
To describe the adiabatic evolution associated with the \(j\)-th instantaneous spectral subspace, we introduce the adiabatic intertwiner generated by \(H_{A,j}(s)\). For each \(j\in \mathbf{I}\), let
\[
U_{A,j}:[0,1]\to \mathcal{U}(\mathcal{H})
\]
denote the solution of the initial-value problem
\begin{equation}
i\,\frac{d}{ds}U_{A,j}(s)=H_{A,j}(s)\,U_{A,j}(s),
\qquad s\in[0,1],
\label{eq:13}
\end{equation}
\[
U_{A,j}(0)=I.
\]
Its adjoint evolution is given by
\begin{equation}
i\,\frac{d}{ds}U_{A,j}(s)^*
=
-\,U_{A,j}(s)^*\,H_{A,j}(s),
\qquad s\in[0,1],
\label{eq:14}
\end{equation}
\[
U_{A,j}(0)^*=I.
\]
\end{definition}
\begin{theorem}[Resolvent-based operator estimates]
\label{thm:resolvent-operator-estimates}
Let \(H_T\) be an operator-valued function satisfying the Adiabatic Assumption and let
\(X_j[A](s)\) be as in \eqref{eq:8}. Then the following bounds hold for any self-adjoint bounded operator \(A\in C^1([0,1])\):
\begin{equation}
\|X_j[A](s)\|
\le
\frac{\|A(s)\|}{g_j(s)},
\qquad s\in[0,1].
\label{eq:15}
\end{equation}
such that, 
\[
X_j[A](s):[0,1]\to \mathcal{B}_{\mathrm{sa}}(\mathcal{H}).
\]
Furthermore,
\begin{equation}
\left\|\frac{d}{ds}X_j[A](s)\right\|
\le
\frac{\|\dot{A}(s)\|}{g_j(s)}
+\frac{4\,\|\dot{H}_T(s)\|\,\|A(s)\|}{g_j(s)^2},
\qquad s\in[0,1].
\label{eq:16}
\end{equation}
\end{theorem}
\begin{proof}
See~\cite[pp.~632--635]{Scherer2019}.
\end{proof}
\begin{lemma}[Norm of the commutator for spectral projection and its derivatives]
\label{lem:spectral-projection-norm}
Let \(H_T\) be an operator-valued function satisfying the Adiabatic Assumption and let $\{P_j(s)\}_{j\in \mathbf{I}}$ denote the family of orthogonal spectral projections onto the eigenspaces of $H_T(s)$. Then,
\begin{equation} 
\|[\dot P_j(s),P_j(s)]\| = \|\dot P_j(s)\|,
\qquad s\in[0,1].
\label{eq:17} 
\end{equation}
\end{lemma}
\begin{proof}
The orthogonal spectral projection \(P_j(s)\) of \(H_T(s)\) satisfies \(P_j(s)^2=P_j(s)\).\\
Differentiating this identity yields
\[
\dot P_j(s)P_j(s)+P_j(s)\dot P_j(s)=\dot P_j(s).
\]
Consequently,
\begin{align*}
[\dot P_j(s),P_j(s)]
&=\dot P_j(s)P_j(s)-P_j(s)\dot P_j(s)\notag\\
\quad
&=2\dot P_j(s)P_j(s)-\dot P_j(s)\notag\\
\quad
&=\dot P_j(s)\bigl(2P_j(s)-I\bigr).
\end{align*}
Since \(P_j(s)\) is an orthogonal spectral projection, the operator
\(2P_j(s)-I\) is self-adjoint.\\
Then,
\begin{align*}
\bigl(2P_j(s)-I\bigr)\bigl(2P_j(s)-I\bigr)^*
&=\bigl(2P_j(s)-I\bigr)^2\notag\\
\quad
&=4P_j(s)^2-4P_j(s)+I\notag\\
\quad
&=4P_j(s)-4P_j(s)+I\notag\\
\quad
&=I.
\end{align*}
Similarly,
\[
\bigl(2P_j(s)-I\bigr)^*\bigl(2P_j(s)-I\bigr)
=I.
\]
Thus \(2P_j(s)-I\) is unitary.

\medskip
\noindent
By invariance of the operator norm
under right multiplication by a unitary operator,
\begin{align*}
\bigl\|[\dot P_j(s),P_j(s)]\bigr\|
&=\bigl\|\dot P_j(s)\bigl(2P_j(s)-I\bigr)\bigr\|\notag\\
\quad
&=\|\dot P_j(s)\|.
\end{align*}
\end{proof}
\begin{corollary}[Bound for spectral projection derivatives]
\label{cor:spectral-projection-derivative-bound}
Let \(H_T\) be an operator-valued function satisfying the Adiabatic Assumption. Let $\{P_j(s)\}_{j\in \mathbf{I}}$ denote the family of orthogonal projections onto the eigenspaces of $H_T(s)$ and $g_j(s)$ be the corresponding spectral gap. Then,
\begin{equation}
\|\dot P_j(s)\| \le \frac{\|\dot H_T(s)\|}{g_j(s)},
\qquad s\in[0,1].
\label{eq:18}
\end{equation}
\end{corollary}
\begin{proof}
See~\cite[Corollary~G.14, p.~636]{Scherer2019}.
\end{proof}
\begin{lemma}[Double commutator identity and norm bound]
\label{lem:double-commutator-norm}
Let \(H_T\) be an operator-valued function satisfying the Adiabatic Assumption and for each $j\in\mathbf{I}$, let $P_j(s)$ be an orthogonal spectral projection onto the eigenspace of $H_T(s)$.
Then
\begin{equation}
\left[H_T(s),\bigl[H_T(s),X_j[\dot{P}_j](s)\bigr]\right]
=
\left[\bigl[\dot H_T(s),P_j(s)\bigr],P_j(s)\right],
\qquad s\in[0,1].
\label{eq:19}
\end{equation}
Furthermore,
\begin{equation}
\left\|
\Bigl[\bigl[\dot H_T(s),P_j(s)\bigr],P_j(s)\Bigr]
\right\|
\leq
\|\dot H_T(s)\|,
\qquad s\in[0,1].
\label{eq:20}
\end{equation}
\end{lemma}
\begin{proof}
\(P_j(s)\) is an orthogonal spectral projection of \(H_T(s)\).\\
Therefore,
\[
[H_T(s),P_j(s)]=0.
\]
Differentiating this identity yields
\[
[H_T(s),\dot P_j(s)]
=-[\dot H_T(s),P_j(s)].
\]
From \eqref{eq:11} and the Jacobi identity, we obtain
\begin{align*}
\left[H_T(s),\bigl[H_T(s),X_j[\dot P_j](s)\bigr]\right]
&=-\left[H_T(s),[\dot P_j(s),P_j(s)]\right]\\
&=-\left[[H_T(s),\dot P_j(s)],P_j(s)\right]\\
&=\left[[\dot H_T(s),P_j(s)],P_j(s)\right].
\end{align*}
The orthogonal spectral projection of \(H_T(s)\) identity \(P_j(s)^2=P_j(s)\) gives
\begin{align*}
\left[[\dot H_T(s),P_j(s)],P_j(s)\right]
&=\dot H_T(s)P_j(s)+P_j(s)\dot H_T(s)
  -2P_j(s)\dot H_T(s)P_j(s)\\
&=\frac12\left(
\dot H_T(s)
-\bigl(2P_j(s)-I\bigr)\dot H_T(s)\bigl(2P_j(s)-I\bigr)
\right).
\end{align*}
\(2P_j(s)-I\) is unitary as shown in the proof of \eqref{eq:17}.

\medskip
\noindent
From the triangle inequality and unitary invariance of the operator norm,
\begin{align*}
\left\|\left[[\dot H_T(s),P_j(s)],P_j(s)\right]\right\|
&\leq
\frac{1}{2}\left(\|\dot H_T(s)\|+\|\dot H_T(s)\|\right)\notag\\
\quad
&\leq
\|\dot H_T(s)\|.
\end{align*}
\end{proof}
\begin{theorem}[Quantum Adiabatic Theorem]
\label{thm:quantum adiabatic}
Let $H_T$ be as defined in \eqref{eq:2}, such that it satisfy the Adiabatic Assumption, and let $\{P_i(s)\}_{i\in \mathbf{I}}$ be the family of orthogonal spectral projections onto the eigenspaces of $H_T(s)$. For a given $j\in \mathbf{I}$, let $g_j$  be its gap function as defined in \eqref{eq:9}. Moreover, let $U_T$ and $U_{A,j}$ be as defined in \eqref{eq:5} and \eqref{eq:13}.
Then, for all $j\in \mathbf{I}$ and $s\in[0,1]$,
\begin{equation}
\left\| \bigl(U_{A,j}(s)-U_T(s)\bigr)P_j(0) \right\|
\leq \frac{C_j(s)}{T},
\label{eq:21}
\end{equation}
where
\begin{equation}
C_j(s)
:=
\frac{\|\dot{H}_T(s)\|}{g_j(s)^2}
+\frac{\|\dot{H}_T(0)\|}{g_j(0)^2}
+\int_0^s\left(\frac{\|\ddot{H}_T(u)\|}{g_j(u)^2}
+10\,\frac{\|\dot{H}_T(u)\|^2}{g_j(u)^3}\right)\,du.
\label{eq:22}
\end{equation}
\end{theorem}
This provides a deterministic bound on the error of the adiabatic approximation for initial states in the eigenspace \(P_j(0)\mathcal H\). In particular, it shows that the difference between the exact propagator $U_T(s)$ and the adiabatic intertwiner $U_{A,j}(s)$ is \(O(T^{-1})\), with a constant that depends explicitly on the first and second derivatives of the Hamiltonian and on inverse powers of the spectral gap.
\begin{proof}
See~\cite[pp.~637--640]{Scherer2019}.
\end{proof}
\section{Propagator Method}
The propagator method studies stochastic adiabatic evolution by comparing the exact stochastic propagator \(U_T(s)\) with the stochastic adiabatic intertwiner \(U_{A,j}(s)\). This comparison provides an operator-level description of the error in the adiabatic approximation. The propagator method is formulated on a finite-dimensional Hilbert space
\(\mathcal{H}\).\\
Introduce the rescaled time \(s\in[0,1]\) by
\begin{align}
t
&=
t_{\mathrm{ini}}+Ts\notag\\
\quad
dt
&=
T\,ds.
\label{eq:23}
\end{align}
Let \((\Omega,\mathcal{F},\{\mathcal{F}_t\}_{t\in[t_{\mathrm{ini}},t_{\mathrm{fin}}]},\mathbb{P})\) be a filtered probability space satisfying the usual conditions, and 
let \(\{W_t\}_{t\in[t_{\mathrm{ini}},t_{\mathrm{fin}}]}\) be a standard real-valued Wiener process adapted to the filtration \(\mathcal{F}_t\).\\
To preserve the Wiener covariance under the change of time, define the rescaled Wiener process by
\begin{align}
W_s
&=
\frac{W_{t_{\mathrm{ini}}+Ts}-W_{t_{\mathrm{ini}}}}{\sqrt{T}}\notag\\
\quad
dW_t
&=
\sqrt{T}\,dW_s.
\label{eq:24}
\end{align}
By the scaling and stationary increment properties of the Wiener process, 
\(\{W_s\}_{s\in[0,1]}\) is again a standard Wiener process with respect to the rescaled filtration 
\[
\mathcal{F}_T(s) = \mathcal{F}_{t_{\mathrm{ini}}+Ts}.
\]
Throughout the stochastic adiabatic analysis, both the exact propagator \(U_T(s)\) 
and the adiabatic intertwiner \(U_{A,j}(s)\) are defined on the same filtered probability space and are driven by the same rescaled Wiener process \(W_s\). Therefore, the stochastic increment \(dW_s\) appearing in the stochastic differential equations for \(U_T(s)\) and \(U_{A,j}(s)\) represents the same Wiener increment.

\bigskip
\noindent
Assume that the stochastic Schrödinger dynamics is written directly in Itô form as
\begin{equation*}
d\psi(t)
=
\left(-iH(t)-\frac{\alpha^2}{2}H(t)^2\right)\psi(t)dt
-i\alpha H(t)\psi(t)dW_t,
\end{equation*}
where  \(\{W_t\}_{t\in[t_{\mathrm{ini}},t_{\mathrm{fin}}]}\) is a standard real-valued Wiener process and \(\alpha\geq0\) denotes the physical noise amplitude.

\medskip
\noindent
Define the state as
\[
\psi(t)=U(t)\psi_{\mathrm{ini}}.
\]
Then the corresponding evolution operator \(U(t)\) satisfies the stochastic differential equation
\begin{equation}
dU(t)
=
A(t)U(t)\,dt
+B(t)U(t)\,dW_t,
\label{eq:25}
\end{equation}
with
\begin{equation}
A(t)
=
-iH(t)
-\frac{\alpha^2}{2}H(t)^2,
\qquad
B(t)
=
-i\alpha H(t).
\label{eq:26}
\end{equation}
Thus, \(A(t)\) is the drift operator and \(B(t)\) is the diffusion operator. In the present model, the diffusion operator is proportional to the instantaneous Hamiltonian,
\begin{equation*}
B(t)\propto H(t).
\end{equation*}
Consequently, the stochastic perturbation is aligned with the instantaneous Hamiltonian at each time.

\bigskip
\noindent
To describe weak noise
over an evolution time \(T\), let
\begin{equation}
\alpha=\sigma T^{-\gamma},
\qquad
\sigma\geq 0,\quad \gamma\geq 0,
\label{eq:27}
\end{equation}
where \(\sigma\) and \(\gamma\) are independent of \(T\).\\
For every \(\gamma>0\),
\[
\alpha\longrightarrow 0
\qquad
\text{as}
\qquad
T\longrightarrow\infty,
\]
so the instantaneous stochastic perturbation becomes weak in physical
time. In contrast, \(\gamma=0\) gives \(\alpha=\sigma\), corresponding to a fixed physical noise amplitude.
So,
\begin{equation}
A(t)=-iH(t)-\frac{\sigma^2T^{-2\gamma}}{2}H(t)^2,
\qquad
B(t)=-i\sigma T^{-\gamma} H(t).
\label{eq:28}
\end{equation}
Define
\begin{equation*}
H_T(s)=H(t_{\mathrm{ini}}+Ts),
\qquad
U_T(s)=U(t_{\mathrm{ini}}+Ts,t_{\mathrm{ini}}).
\end{equation*}
The stochastic differential equation \eqref{eq:25} in rescaled time becomes
\begin{equation}  
dU_T(s)=T\,A_T(s)\,U_T(s)\,ds+\sqrt{T}\,B_T(s)\,U_T(s)\,dW_s.
\label{eq:29}
\end{equation}
Here,
\begin{equation}
A_T(s)=-iH_T(s)-\frac{\sigma^2T^{-2\gamma}}{2}H_T(s)^2,
\qquad
B_T(s)=-i\sigma T^{-\gamma}H_T(s).
\label{eq:30}
\end{equation}
Consequently, the exact evolution \(U_T(s)\) is governed by the stochastic differential equation
\begin{align}
dU_T(s)
=
&\left(-iT H_T(s)-\frac{\sigma^2T^{1-2\gamma}}{2}H_T(s)^2\right)U_T(s)\,ds\notag\\
\quad
&-i\sigma T^{\frac{1}{2}-\gamma}H_T(s)U_T(s)\,dW_s,
\qquad
U_T(0)=I.
\label{eq:31}
\end{align}
Substituting \eqref{eq:12} into \eqref{eq:13} gives
\begin{align}
dU_{A,j}(s)
&=
-i\left(T H_T(s)+i\bigl[\dot{P}_j(s),\,P_j(s)\bigr]\right)U_{A,j}(s)\,ds\notag\\
\quad
&=
\left(-iT H_T(s)+\bigl[\dot{P}_j(s),\,P_j(s)\bigr]\right)U_{A,j}(s)\,ds.
\label{eq:32}
\end{align}
Accordingly, we define the stochastic adiabatic intertwiner \(U_{A,j}(s)\) directly through the stochastic differential equation
\begin{align}
dU_{A,j}(s)
&=
\left(-iT H_T(s)+\bigl[\dot{P}_j(s),\,P_j(s)\bigr]
-\frac{\sigma^2T^{1-2\gamma}}{2}H_T(s)^2\right)U_{A,j}(s)\,ds\notag\\
&\quad
-i\sigma T^{\frac{1}{2}-\gamma}H_T(s)U_{A,j}(s)\,dW_s.
\label{eq:33}
\end{align}
Using \eqref{eq:12}, this becomes
\begin{align}
dU_{A,j}(s)
&=
\left[-iH_{A,j}(s)-\frac{\sigma^2T^{1-2\gamma}}{2}H_T(s)^2\right]U_{A,j}(s)\,ds\notag\\
&\quad
-i\sigma T^{\frac{1}{2}-\gamma}H_T(s)U_{A,j}(s)\,dW_s,
\qquad 
U_{A,j}(0)=I.
\label{eq:34}
\end{align}
Then the equation for the stochastic adiabatic intertwiner can be further reduced to
\begin{equation}
dU_{A,j}(s)
=
\left(TA_T(s)+\bigl[\dot{P}_j(s),\,P_j(s)\bigr]\right)U_{A,j}(s)\,ds
+\sqrt{T}\,B_T(s)U_{A,j}(s)\,dW_s.
\label{eq:35}
\end{equation}
Thus \(U_T(s)\) and \(U_{A,j}(s)\) contain the same stochastic dynamical terms, while their drift generators differ only by the term \(\bigl[\dot{P}_j(s),\,P_j(s)\bigr]\).
\begin{lemma}
\label{lem:stochastic-intertwining}
Let \(U_{A,j}(s)\) be the stochastic adiabatic intertwiner defined by
\eqref{eq:33}, and let \(P_j(s)\) be the orthogonal spectral projection associated
with the \(j\)-th eigenspace of \(H_T(s)\). Suppose that \(H_T(s)\) satisfies the Adiabatic
Assumption and
\[
U_{A,j}(0)=I.
\]
Then \(U_{A,j}(s)\) intertwines the initial and instantaneous spectral projections.
\begin{equation}
P_j(s)U_{A,j}(s)
=
U_{A,j}(s)P_j(0),
\qquad s\in[0,1]
\qquad \text{a.s..}
\label{eq:36}
\end{equation}
\end{lemma}
\begin{proof}
Let \(P_j(s)\) be the orthogonal spectral projection associated with the self-adjoint Hamiltonian \(H_T(s)\).\\
Then,
\begin{equation*}
[P_j(s),H_T(s)]=0.
\end{equation*}
Hence \(P_j(s)\) also commutes with \(H_T(s)^2\).
Therefore,
\begin{equation}
[P_j(s),A_T(s)]=0,
\qquad
[P_j(s),B_T(s)]=0.
\label{eq:37}
\end{equation}
Since \(P_j(s)^2=P_j(s)\), differentiation gives
\begin{equation}
\dot P_j(s)P_j(s)+P_j(s)\dot P_j(s)=\dot P_j(s).
\label{eq:38}
\end{equation}
Multiplying by \(P_j(s)\) on both sides yields
\begin{equation*}
P_j(s)\dot P_j(s)P_j(s)=0.
\end{equation*}
Consequently,
\begin{align*}
P_j(s)\bigl[\dot P_j(s),\,P_j(s)\bigr]
&=
P_j(s)\dot P_j(s)P_j(s)-P_j(s)\dot P_j(s)\notag\\
\quad
&=
-P_j(s)\dot P_j(s),
\end{align*}
whereas
\begin{equation*}
\bigl[\dot P_j(s),\,P_j(s)\bigr]P_j(s)=\dot P_j(s)P_j(s).
\end{equation*}
Combining these results with \eqref{eq:38}, we obtain
\begin{equation}
\dot P_j(s)+P_j(s)\bigl[\dot P_j(s),\,P_j(s)\bigr]
=\bigl[\dot P_j(s),\,P_j(s)\bigr]P_j(s).
\label{eq:39}
\end{equation}
Define the intertwining error process
\begin{equation}
F_j(s)
=
P_j(s)U_{A,j}(s)-U_{A,j}(s)P_j(0).
\label{eq:40}
\end{equation}
This immediately gives
\begin{equation*}
F_j(0)=0.
\end{equation*}
Apply It\^o's product formula to \(F_j(s)\)
\begin{equation*}
dF_j(s)
=
dP_j(s)U_{A,j}(s)
+P_j(s)\,dU_{A,j}(s)
+dP_j(s)\,dU_{A,j}(s)
-dU_{A,j}(s)\,P_j(0).
\end{equation*}
Because \(P_j(s)\) is deterministic and of finite variation,
\[
dP_j(s)\,dU_{A,j}(s)=0.
\]
Thus,
\begin{align}
dF_j(s)
=
\dot P_j(s)U_{A,j}(s)\,ds
+P_j(s)\,dU_{A,j}(s)
-\,dU_{A,j}(s)P_j(0).
\label{eq:41}
\end{align}
Substituting \eqref{eq:35} and \eqref{eq:40}, we obtain
\begin{align*}
dF_j(s)
=
&\,TA_T(s)F_j(s)\,ds
+\Bigl(\dot P_j(s)+P_j(s)\bigl[\dot P_j(s),\,P_j(s)\bigr]\Bigr)U_{A,j}(s)\,ds\notag\\
\quad
&-\bigl[\dot P_j(s),\,P_j(s)\bigr]U_{A,j}(s)P_j(0)\,ds
+\sqrt{T}\,B_T(s)F_j(s)\,dW_s.
\end{align*}
From \eqref{eq:39},
\begin{align*}
dF_j(s)
=
&\,TA_T(s)F_j(s)\,ds
+\bigl[\dot P_j(s),\,P_j(s)\bigr]
\Bigl(P_j(s)U_{A,j}(s)-U_{A,j}(s)P_j(0)\Bigr)\,ds\notag\\
\quad
&+\sqrt{T}\,B_T(s)F_j(s)\,dW_s.
\end{align*}
Therefore, by \eqref{eq:40},
\begin{equation*}
dF_j(s)
=
\left(TA_T(s)+\bigl[\dot P_j(s),\,P_j(s)\bigr]\right)F_j(s)\,ds
+\sqrt{T}\,B_T(s)F_j(s)\,dW_s.
\end{equation*}
This is a homogeneous linear It\^o stochastic differential equation.
The zero process is a solution, and uniqueness of the solution with
initial condition \(F_j(0)=0\) implies
\[
F_j(s)=0
\]
for every \(s\in[0,1]\), \text{a.s..}\\
Hence
\begin{equation*}
P_j(s)U_{A,j}(s)
=
U_{A,j}(s)P_j(0),
\qquad s\in[0,1].
\qquad \text{a.s..}
\end{equation*}

\medskip
\noindent
Consequently, if
\[
\psi_0\in\operatorname{Ran}P_j(0),
\]
then
\begin{equation*}
U_{A,j}(s)\psi_0
\in
\operatorname{Ran}P_j(s)
\qquad s\in[0,1].
\qquad \text{a.s..}
\end{equation*}
\end{proof}
\begin{lemma}
\label{lem:almost-sure-unitarity}
Let \(U_T(s)\) and \(U_{A,j}(s)\) denote the exact stochastic propagator
and the stochastic adiabatic intertwiner, defined by \eqref{eq:31} and
\eqref{eq:33}, respectively, with
\[
U_T(0)=U_{A,j}(0)=I.
\]
Then, almost surely, \(U_T(s)\) and \(U_{A,j}(s)\) are unitary for every
\(s\in[0,1]\). Equivalently,
\begin{align}
U_T(s)^*U_T(s)
&=U_T(s)U_T(s)^*=I,\notag\\
\quad
U_{A,j}(s)^*U_{A,j}(s)
&=U_{A,j}(s)U_{A,j}(s)^*=I,
\qquad s\in[0,1]
\qquad \text{a.s..}
\label{eq:42}
\end{align}
\end{lemma}
\begin{proof}
Let \(P_j(s)\) be an orthogonal spectral projection of the self-adjoint \(H_T(s)\).\\
Taking the adjoint of \eqref{eq:29} gives
\begin{equation}
dU_T(s)^*
=
T\,U_T(s)^*A_T(s)^*\,ds
+
\sqrt{T}\,U_T(s)^*B_T(s)^*\,dW_s.
\label{eq:43}
\end{equation}
Applying It\^o's product formula to \(U_T(s)^*U_T(s)\),
\begin{align}
d\left(U_T(s)^*U_T(s)\right)
&=
dU_T(s)^*\,U_T(s)
+U_T(s)^*\,dU_T(s)
+dU_T(s)^*\,dU_T(s).
\label{eq:44}
\end{align}
The quadratic-variation term is
\begin{equation}
dU_T(s)^*\,dU_T(s)
=
T\,U_T(s)^*B_T(s)^*B_T(s)U_T(s)\,ds.
\label{eq:45}
\end{equation}
From \eqref{eq:29}, \eqref{eq:43} and \eqref{eq:45},
\begin{align}
d\left(U_T(s)^*U_T(s)\right)
&=
T\,U_T(s)^*\left(A_T(s)+A_T(s)^*+B_T(s)^*B_T(s)\right)U_T(s)\,ds\notag\\
&\quad
+\sqrt{T}\,U_T(s)^*\left(B_T(s)^*+B_T(s)\right)U_T(s)\,dW_s.
\label{eq:46}
\end{align}
Since \(H_T(s)\) is self-adjoint,
\begin{align*}
A_T(s)^*
&=
iH_T(s)-\frac{\sigma^2T^{-2\gamma}}{2}H_T(s)^2,\\
A_T(s)+A_T(s)^*
&=
-\sigma^2T^{-2\gamma}H_T(s)^2,\\
B_T(s)^*
&=
i\sigma T^{-\gamma}H_T(s),\\
B_T(s)^*B_T(s)
&=
\sigma^2T^{-2\gamma}H_T(s)^2.
\end{align*}
Therefore,
\begin{equation}
A_T(s)
+A_T(s)^*
+B_T(s)^*B_T(s)
=
0.
\label{eq:47}
\end{equation}
Moreover,
\begin{equation}
B_T(s)^*+B_T(s)=0.
\label{eq:48}
\end{equation}
Substituting \eqref{eq:47} and \eqref{eq:48} into \eqref{eq:46} yields
\begin{equation*}
d\left(U_T(s)^*U_T(s)\right)=0.
\end{equation*}
Thus \(U_T(s)^*U_T(s)\) is constant in \(s\in[0,1]\).

\medskip
\noindent
Since \(U_T(0)=I\),
\begin{equation*}
U_T(s)^*U_T(s)=I
\qquad s\in[0,1].
\qquad \text{a.s..}
\end{equation*}

\bigskip
\noindent
We next consider the stochastic adiabatic intertwiner \(U_{A,j}(s)\).
Since \(P_j(s)\) is an orthogonal projection,
\[
P_j(s)^*=P_j(s),
\qquad
\dot P_j(s)^*=\dot P_j(s).
\]
Therefore,
\begin{align*}
\bigl[\dot P_j(s),P_j(s)\bigr]^*
&=
\left(\dot P_j(s)P_j(s)-P_j(s)\dot P_j(s)\right)^*\notag\\
\quad
&=
P_j(s)\dot P_j(s)-\dot P_j(s)P_j(s)\notag\\
\quad
&=
-\bigl[\dot P_j(s),P_j(s)\bigr].
\end{align*}
Hence
\begin{equation}
\bigl[\dot P_j(s),P_j(s)\bigr]^*+\bigl[\dot P_j(s),P_j(s)\bigr]=0.
\label{eq:49}
\end{equation}
Taking adjoints of \eqref{eq:35} gives
\begin{equation}
dU_{A,j}(s)^*
=
U_{A,j}(s)^*\left(T\,A_T(s)^*+\bigl[\dot P_j(s),P_j(s)\bigr]^*\right)ds
+\sqrt{T}\,U_{A,j}(s)^*B_T(s)^*\,dW_s.
\label{eq:50}
\end{equation}
Apply It\^o's product formula to \(U_{A,j}(s)^*U_{A,j}(s)\),
\begin{align}
d\left(U_{A,j}(s)^*U_{A,j}(s)\right)
&=
dU_{A,j}(s)^*\,U_{A,j}(s)
+U_{A,j}(s)^*\,dU_{A,j}(s)
+dU_{A,j}(s)^*\,dU_{A,j}(s).
\label{eq:51}
\end{align}
The quadratic-variation term is
\begin{equation}
dU_{A,j}(s)^*\,dU_{A,j}(s)
=
T\,U_{A,j}(s)^*B_T(s)^*B_T(s)U_{A,j}(s)\,ds.
\label{eq:52}
\end{equation}
Obtain by \eqref{eq:35}, \eqref{eq:50} and \eqref{eq:52},
\begin{align}
d\left(U_{A,j}(s)^*U_{A,j}(s)\right)
&=
U_{A,j}(s)^*\Bigl[T\bigl(A_T(s)+A_T(s)^*+B_T(s)^*B_T(s)\bigr)\notag\\
&\quad
+\bigl[\dot P_j(s),P_j(s)\bigr]^*+\bigl[\dot P_j(s),P_j(s)\bigr]\Bigr]
U_{A,j}(s)\,ds\notag\\
&\quad
+\sqrt{T}\,U_{A,j}(s)^*\bigl(B_T(s)^*+B_T(s)\bigr)U_{A,j}(s)\,dW_s.
\label{eq:53}
\end{align}
From \eqref{eq:47},\eqref{eq:48}, and \eqref{eq:49}, this becomes,
\begin{equation*}
d\left(U_{A,j}(s)^*U_{A,j}(s)\right)=0.
\end{equation*}
Consequently, \(U_{A,j}(s)^*U_{A,j}(s)\) is constant in \(s\in[0,1]\).

\medskip
\noindent
Since \(U_{A,j}(0)=I\), it follows that
\begin{equation*}
U_{A,j}(s)^*U_{A,j}(s)=I,
\qquad 
s\in[0,1],
\end{equation*}
almost surely.

\medskip
\noindent
Finally, because the Hilbert space \(\mathcal H\) is finite-dimensional, each isometry
on \(\mathcal H\) is unitary. Therefore,
\[
U_T(s)U_T(s)^*
=
U_T(s)^*U_T(s)
=
I
\]
and
\[
U_{A,j}(s)U_{A,j}(s)^*
=
U_{A,j}(s)^*U_{A,j}(s)
=
I
\qquad s\in[0,1].
\qquad \text{a.s..}
\]
\end{proof}
\begin{theorem}[Stochastic Version of Quantum Adiabatic Theorem]
\label{thm:stochastic version quantum adiabatic}
Let $H_T$ be as defined in \eqref{eq:2}, acting on a Hilbert space $\mathcal H$ of fixed finite dimension, and satisfying the Adiabatic Assumption. Let $\{P_i(s)\}_{i\in \mathbf{I}}$ be the family of orthogonal projections onto the eigenspaces of $H_T(s)$. For a given $j\in \mathbf{I}$, let $g_j$ be its gap function as defined in \eqref{eq:9}.

\medskip
\noindent
Let
\[
\alpha=\sigma T^{-\gamma},
\qquad \sigma\geq 0,\quad \gamma\geq 0,
\]
be the physical noise amplitude in \eqref{eq:27}, where
$\sigma$ and $\gamma$ are independent of $T$. Moreover, let $U_T$ and $U_{A,j}$ be the exact stochastic propagator and stochastic adiabatic intertwiner defined by \eqref{eq:31} and \eqref{eq:33}. Both operators are driven by the same standard real-valued Wiener process on a filtered probability space satisfying the usual conditions, and $U_T(0)=U_{A,j}(0)=I$.

\medskip
\noindent
Then, for every $T\geq1$,
\begin{align}
\E\sup_{0\le s\le 1}\|(U_{A,j}(s)-U_T(s))P_j(0)\|^2
&\le
\frac{3}{T^2}\sup_{0\le s\le 1}C_j(s)^2
+\frac{3\sigma^4\,T^{-4\gamma}}{4}
\Bigl(\int_0^1\|\dot H_T(u)\|\,du\Bigr)^2\notag\\
&\quad
+12\sigma^2\,d\,T^{-(2\gamma+1)}
\int_0^1\frac{\|\dot H_T(u)\|^2}{g_j(u)^2}\ du
\label{eq:54}
\end{align}
where
\begin{equation}
C_j(s)
=
\frac{\|\dot H_T(s)\|}{g_j(s)^2}
+\frac{\|\dot H_T(0)\|}{g_j(0)^2}
+\int_0^s
\left(\frac{\|\ddot H_T(u)\|}{g_j(u)^2}+10\,\frac{\|\dot H_T(u)\|^2}{g_j(u)^3}
\right)\,du.
\label{eq:55}
\end{equation}
\end{theorem}
\begin{proof}
Apply It\^o's product formula to \(U_{A,j}(s)^*U_T(s)\),
\begin{equation}
d\bigl(U_{A,j}(s)^*U_T(s)\bigr)
=
dU_{A,j}(s)^*U_T(s)
+U_{A,j}(s)^*\,dU_T(s)
+dU_{A,j}(s)^*dU_T(s).
\label{eq:56}
\end{equation}
From \eqref{eq:29} and \eqref{eq:35}, we obtain,
\begin{align}
dU_{A,j}(s)^*dU_T(s)
&=
\left(\sqrt{T}\,U_{A,j}(s)^*B_T(s)^*\,dW_s\right)
\left(\sqrt{T}\,B_T(s)U_T(s)\,dW_s\right)\notag\\
\quad
&=
T\,U_{A,j}(s)^*B_T(s)^*B_T(s)U_T(s)\,ds.
\label{eq:57}
\end{align}
Substituting \eqref{eq:29}, \eqref{eq:50} and \eqref{eq:57} into \eqref{eq:56} gives
\begin{align*}
d\bigl(U_{A,j}(s)^*U_T(s)\bigr)
&=
U_{A,j}(s)^*
\left[TA_T(s)^*-\bigl[\dot{P}_j(s),\,P_j(s)\bigr]\right]U_T(s)\,ds\notag\\
&\quad
+\sqrt{T}\,U_{A,j}(s)^*B_T(s)^*U_T(s)\,dW_s
+T\,U_{A,j}(s)^*A_T(s)U_T(s)\,ds\notag\\
&\quad
+\sqrt{T}\,U_{A,j}(s)^*B_T(s)U_T(s)\,dW_s
+T\,U_{A,j}(s)^*B_T(s)^*B_T(s)U_T(s)\,ds.
\end{align*}
Factoring \(U_{A,j}(s)^*\) from the left and \(U_T(s)\) from the right yields
\begin{align}
d\bigl(U_{A,j}(s)^*U_T(s)\bigr)
&=
U_{A,j}(s)^*\left[T\bigl(A_T(s)^*+A_T(s)+B_T(s)^*B_T(s)\bigr)
-\bigl[\dot{P}_j(s),\,P_j(s)\bigr]\right]U_T(s)\,ds\notag\\
&\quad
+
\sqrt{T}\,U_{A,j}(s)^*\bigl(B_T(s)^*+B_T(s)\bigr)U_T(s)\,dW_s.
\label{eq:58}
\end{align}
Obtain by \eqref{eq:47} and \eqref{eq:48}
\begin{equation*}
d\bigl(U_{A,j}(s)^*U_T(s)\bigr)
=
-U_{A,j}(s)^*\bigl[\dot P_j(s),P_j(s)\bigr]U_T(s)\,ds.
\end{equation*}
From \eqref{eq:11},
\begin{equation*}
d\bigl(U_{A,j}(s)^*U_T(s)\bigr)
=
U_{A,j}(s)^*[H_T(s),X_j[\dot P_j](s)]U_T(s)\,ds.
\end{equation*}
Expanding the commutator gives
\begin{align*}
d\left(U_{A,j}(s)^*U_T(s)\right)
&=
U_{A,j}(s)^*\bigl(H_T(s)X_j[\dot P_j](s)-X_j[\dot P_j](s)H_T(s)\bigr)U_T(s)\,ds.
\end{align*}
By \eqref{eq:12},
\begin{align*}
d\!\left(U_{A,j}(s)^*U_T(s)\right)
&=
U_{A,j}(s)^*
\left[\frac{1}{T}\Bigl(H_{A,j}(s)-i\bigl[\dot P_j(s),P_j(s)\bigr]\Bigr)
X_j[\dot P_j](s)-X_j[\dot P_j](s)H_T(s)\right]U_T(s)\,ds.
\end{align*}
It follows that
\begin{align}
d\!\left(U_{A,j}(s)^*U_T(s)\right)
&=
\frac{1}{T}U_{A,j}(s)^*H_{A,j}(s)X_j[\dot P_j](s)U_T(s)\,ds
-\frac{i}{T}U_{A,j}(s)^*\bigl[\dot P_j(s),P_j(s)\bigr]X_j[\dot P_j](s)U_T(s)\,ds\notag\\
&\quad
-U_{A,j}(s)^*X_j[\dot P_j](s)H_T(s)U_T(s)\,ds.
\label{eq:59}
\end{align}
Substituting \eqref{eq:31} and \eqref{eq:34} into \eqref{eq:59}, we obtain
\begin{align*}
d\!\left(U_{A,j}(s)^*U_T(s)\right)
=
&-\frac{i}{T}\Bigl[dU_{A,j}(s)^*X_j[\dot P_j](s)U_T(s)
+U_{A,j}(s)^*\bigl[\dot P_j(s),P_j(s)\bigr]X_j[\dot P_j](s)U_T(s)\,ds\notag\\
\quad
&+U_{A,j}(s)^*X_j[\dot P_j](s)\,dU_T(s)\Bigr]\notag\\
\quad
&-\frac{i\sigma^2T^{-2\gamma}}{2}U_{A,j}(s)^*\bigl(H_T(s)^2X_j[\dot P_j](s)
+X_j[\dot P_j](s)H_T(s)^2\bigr)U_T(s)\,ds\notag\\
\quad
&+\sigma T^{-(\gamma+\frac{1}{2})}U_{A,j}(s)^*\bigl(X_j[\dot P_j](s)H_T(s)
-H_T(s)X_j[\dot P_j](s)\bigr)U_T(s)\,dW_s.
\end{align*}
Thus,
\begin{align}
d\!\left(U_{A,j}(s)^*U_T(s)\right)
=
&-\frac{i}{T}\Bigl[dU_{A,j}(s)^*X_j[\dot P_j](s)U_T(s)
+U_{A,j}(s)^*\bigl[\dot P_j(s),P_j(s)\bigr]X_j[\dot P_j](s)U_T(s)\,ds\notag\\
\quad
&+U_{A,j}(s)^*X_j[\dot P_j](s)\,dU_T(s)\Bigr]\notag\\
\quad
&-\frac{i\sigma^2T^{-2\gamma}}{2}U_{A,j}(s)^*\bigl(H_T(s)^2X_j[\dot P_j](s)
+X_j[\dot P_j](s)H_T(s)^2\bigr)U_T(s)\,ds\notag\\
\quad
&-\sigma\,T^{-(\gamma+\frac{1}{2})}U_{A,j}(s)^*\bigl[H_T(s),X_j[\dot P_j](s)\bigr]
U_T(s)\,dW_s.
\label{eq:60}
\end{align}
Now apply It\^o's product formula to \(U_{A,j}(s)^*X_j[\dot P_j](s)U_T(s)\).

\medskip
\noindent
The process \(X_j[\dot P_j](s)\) has finite variation. Therefore, all quadratic products involving \(dX_j[\dot P_j](s)\) vanish.

\medskip
\noindent
Consequently,
\begin{align}
d\left(U_{A,j}(s)^*X_j[\dot P_j](s)U_T(s)\right)
&=
dU_{A,j}(s)^*X_j[\dot P_j](s)U_T(s)
+U_{A,j}(s)^*dX_j[\dot P_j](s)U_T(s)\notag\\
&\quad
+U_{A,j}(s)^*X_j[\dot P_j](s)\,dU_T(s)
+dU_{A,j}(s)^*X_j[\dot P_j](s)dU_T(s).
\label{eq:61}
\end{align}
Obtain by \eqref{eq:29} and \eqref{eq:35},
\begin{align}
dU_{A,j}(s)^*X_j[\dot P_j](s)dU_T(s)
&=
\left(i\sigma\,T^{\frac{1}{2}-\gamma}U_{A,j}(s)^*H_T(s)\,dW_s\right)
X_j[\dot P_j](s)
\left(-i\sigma\,T^{\frac{1}{2}-\gamma}H_T(s)U_T(s)\,dW_s\right)\notag\\
\quad
&=
\sigma^2\,T^{1-2\gamma}U_{A,j}(s)^*H_T(s)X_j[\dot P_j](s)H_T(s)U_T(s)\,ds.
\label{eq:62}
\end{align}
Hence, \eqref{eq:61} becomes
\begin{align}
d\left(U_{A,j}(s)^*X_j[\dot P_j](s)U_T(s)\right)
&=
dU_{A,j}(s)^*X_j[\dot P_j](s)U_T(s)
+U_{A,j}(s)^*dX_j[\dot P_j](s)U_T(s)\notag\\
&\quad
+U_{A,j}(s)^*X_j[\dot P_j](s)\,dU_T(s)\notag\\
&\quad
+\sigma^2\,T^{1-2\gamma}U_{A,j}(s)^*H_T(s)X_j[\dot P_j](s)H_T(s)U_T(s)\,ds.
\label{eq:63}
\end{align}
Substituting \eqref{eq:63} into \eqref{eq:60} yields
\begin{align*}
d\!\left(U_{A,j}(s)^*U_T(s)\right)
=
&-\frac{i}{T}\Bigl[d\left(U_{A,j}(s)^*X_j[\dot P_j](s)U_T(s)\right)
-U_{A,j}(s)^*\frac{d}{ds}X_j[\dot P_j](s)U_T(s)\,ds\notag\\
\quad
&-\sigma^2\,T^{1-2\gamma}U_{A,j}(s)^*H_T(s)X_j[\dot P_j](s)H_T(s)U_T(s)\,ds\notag\\
\quad
&+U_{A,j}(s)^*\bigl[\dot P_j(s),P_j(s)\bigr]X_j[\dot P_j](s)U_T(s)\,ds\Bigr]\notag\\
\quad
&-\frac{i\sigma^2\,T^{-2\gamma}}{2}
U_{A,j}(s)^*\bigr(H_T(s)^2X_j[\dot P_j](s)
+X_j[\dot P_j](s)H_T(s)^2\bigr)U_T(s)\,ds\notag\\
\quad
&-\sigma\,T^{-(\gamma+\frac{1}{2})}U_{A,j}(s)^*\bigl[H_T(s),X_j[\dot P_j](s)\bigr]
U_T(s)\,dW_s.
\end{align*}
Thus,
\begin{align*}
d\!\left(U_{A,j}(s)^*U_T(s)\right)
=
&-\frac{i}{T}\,d\left(U_{A,j}(s)^*X_j[\dot P_j](s)U_T(s)\right)\notag\\
\quad
&-\frac{i}{T}U_{A,j}(s)^*\Bigl(\bigl[\dot P_j(s),P_j(s)\bigr]X_j[\dot P_j](s)
-\frac{d}{ds}X_j[\dot P_j](s)\Bigr)U_T(s)\,ds\notag\\
\quad
&+i\sigma^2\,T^{-2\gamma}U_{A,j}(s)^*H_T(s)X_j[\dot P_j](s)H_T(s)U_T(s)\,ds\notag\\
\quad
&-\frac{i\sigma^2\,T^{-2\gamma}}{2}U_{A,j}(s)^*\bigl(H_T(s)^2X_j[\dot P_j](s)
+X_j[\dot P_j](s)H_T(s)^2\bigr)U_T(s)\,ds\notag\\
\quad
&-\sigma\,T^{-(\gamma+\frac{1}{2})}U_{A,j}(s)^*\bigl[H_T(s),X_j[\dot P_j](s)\bigr]
U_T(s)\,dW_s.
\end{align*}
This gives us
\begin{align*}
d\!\left(U_{A,j}(s)^*U_T(s)\right)
=
&-\frac{i}{T}\,d\left(U_{A,j}(s)^*X_j[\dot P_j](s)U_T(s)\right)\notag\\
\quad
&-\frac{i}{T}U_{A,j}(s)^*\Bigl(\bigl[\dot P_j(s),P_j(s)\bigr]X_j[\dot P_j](s)
-\frac{d}{ds}X_j[\dot P_j](s)\Bigr)U_T(s)\,ds\notag\\
\quad
&+i\sigma^2\,T^{-2\gamma}U_{A,j}(s)^*H_T(s)X_j[\dot P_j](s)H_T(s)U_T(s)\,ds\notag\\
\quad
&-\frac{i\sigma^2\,T^{-2\gamma}}{2}
U_{A,j}(s)^*
\left(
\Bigl[H_T(s),\bigl[H_T(s),X_j[\dot P_j](s)\bigr]\Bigr]\right.\notag\\
\quad
&\left.
+2H_T(s)X_j[\dot P_j](s)H_T(s)\right)U_T(s)\,ds\notag\\
\quad
&-\sigma\,T^{-(\gamma+\frac{1}{2})}U_{A,j}(s)^*\bigl[H_T(s),X_j[\dot P_j](s)\bigr]
U_T(s)\,dW_s.
\end{align*}
Therefore,
\begin{align*}
d\!\left(U_{A,j}(s)^*U_T(s)\right)
=
&-\frac{i}{T}\,d\left(U_{A,j}(s)^*X_j[\dot P_j](s)U_T(s)\right)\notag\\
\quad
&-\frac{i}{T}U_{A,j}(s)^*\Bigl(\bigl[\dot P_j(s),P_j(s)\bigr]X_j[\dot P_j](s)
-\frac{d}{ds}X_j[\dot P_j](s)\Bigr)U_T(s)\,ds\notag\\
\quad
&-\frac{i\sigma^2\,T^{-2\gamma}}{2}U_{A,j}(s)^*\Bigl[H_T(s),\bigl[H_T(s),X_j[\dot P_j](s)\bigr]\Bigr]
U_T(s)\,ds\notag\\
\quad
&-\sigma\,T^{-(\gamma+\frac{1}{2})}U_{A,j}(s)^*\bigl[H_T(s),X_j[\dot P_j](s)\bigr]
U_T(s)\,dW_s.
\end{align*}
By \eqref{eq:11},
\begin{align*}
d\!\left(U_{A,j}(s)^*U_T(s)\right)
=
&-\frac{i}{T}\,d\left(U_{A,j}(s)^*X_j[\dot P_j](s)U_T(s)\right)\notag\\
\quad
&-\frac{i}{T}U_{A,j}(s)^*\Bigl(\bigl[\dot P_j(s),P_j(s)\bigr]X_j[\dot P_j](s)
-\frac{d}{ds}X_j[\dot P_j](s)\Bigr)U_T(s)\,ds\notag\\
\quad
&-\frac{i\sigma^2\,T^{-2\gamma}}{2}U_{A,j}(s)^*
\Bigl[\bigl[\dot H_T(s),P_j(s)\bigr],P_j(s)\Bigr]
U_T(s)\,ds\notag\\
\quad
&+\sigma\,T^{-(\gamma+\frac{1}{2})}U_{A,j}(s)^*\bigl[\dot P_j(s),P_j(s)\bigr]
U_T(s)\,dW_s.
\end{align*}
From \eqref{eq:19},
\begin{align}
d\!\left(U_{A,j}(s)^*U_T(s)\right)
=
&-\frac{i}{T}\,d\left(U_{A,j}(s)^*X_j[\dot P_j](s)U_T(s)\right)\notag\\
\quad
&-\frac{i}{T}U_{A,j}(s)^*\Bigl(\bigl[\dot P_j(s),P_j(s)\bigr]X_j[\dot P_j](s)
-\frac{d}{ds}X_j[\dot P_j](s)\Bigr)U_T(s)\,ds\notag\\
\quad
&-\frac{i\sigma^2\,T^{-2\gamma}}{2}U_{A,j}(s)^*
\Bigl[\bigl[\dot H_T(s),P_j(s)\bigr],P_j(s)\Bigr]
U_T(s)\,ds\notag\\
\quad
&+\sigma\,T^{-(\gamma+\frac{1}{2})}U_{A,j}(s)^*\bigl[\dot P_j(s),P_j(s)\bigr]
U_T(s)\,dW_s.
\label{eq:64}
\end{align}
It follows that
\begin{align*}
U_{A,j}(s)^*U_T(s)
&=
U_{A,j}(0)^*U_T(0)-\frac{i}{T}\Bigl(U_{A,j}(s)^*X_j[\dot P_j](s)U_T(s)
-U_{A,j}(0)^*X_j[\dot P_j](0)U_T(0)\Bigr)\notag\\
&\quad
-\frac{i}{T}\int_0^sU_{A,j}(u)^*\Bigl(\bigl[\dot P_j(u),P_j(u)\bigr]X_j[\dot P_j](u)
-\frac{d}{du}X_j[\dot P_j](u)\Bigr)U_T(u)\,du\notag\\
&\quad
-\frac{i\sigma^2\,T^{-2\gamma}}{2}\int_0^sU_{A,j}(u)^*
\Bigl[\bigl[\dot H_T(u),P_j(u)\bigr],P_j(u)\Bigr]U_T(u)\,du\notag\\
&\quad
+\sigma\,T^{-(\gamma+\frac{1}{2})}\int_0^sU_{A,j}(u)^*\bigl[\dot P_j(u),P_j(u)\bigr]
U_T(u)\,dW_u.
\end{align*} 
Using the initial conditions
\[
U_{A,j}(0)^*=I,
\qquad
U_T(0)=I
\]
we obtain,
\begin{align*}
U_{A,j}(s)^*U_T(s)
&=
I
-\frac{i}{T}\Bigl(U_{A,j}(s)^*X_j[\dot P_j](s)U_T(s)
-X_j[\dot P_j](0)\Bigr)\notag\\
&\quad
-\frac{i}{T}\int_0^sU_{A,j}(u)^*
\Bigl(\bigl[\dot P_j(u),P_j(u)\bigr]X_j[\dot P_j](u)
-\frac{d}{du}X_j[\dot P_j](u)\Bigr)U_T(u)\,du\notag\\
&\quad
-\frac{i\sigma^2\,T^{-2\gamma}}{2}\int_0^sU_{A,j}(u)^*
\Bigl[\bigl[\dot H_T(u),P_j(u)\bigr],P_j(u)\Bigr]U_T(u)\,du\notag\\
&\quad
+\sigma\,T^{-(\gamma+\frac{1}{2})}\int_0^sU_{A,j}(u)^*
\bigl[\dot P_j(u),P_j(u)\bigr]U_T(u)\,dW_u.
\end{align*}
Multiplying from the left by the adiabatic intertwiner $U_{A,j}(s)$ yields
\begin{align*}
U_T(s)
&=
U_{A,j}(s)
-\frac{i}{T}\Bigl(X_j[\dot P_j](s)U_T(s)-U_{A,j}(s)X_j[\dot P_j](0)\Bigr)\notag\\
&\quad
-\frac{i}{T}U_{A,j}(s)\int_0^sU_{A,j}(u)^*
\Bigl(\bigl[\dot P_j(u),P_j(u)\bigr]X_j[\dot P_j](u)
-\frac{d}{du}X_j[\dot P_j](u)\Bigr)U_T(u)\,du\notag\\
&\quad
-\frac{i\sigma^2\,T^{-2\gamma}}{2}U_{A,j}(s)\int_0^sU_{A,j}(u)^*
\Bigl[\bigl[\dot H_T(u),P_j(u)\bigr],P_j(u)\Bigr]U_T(u)\,du\notag\\
&\quad
+\sigma\,T^{-(\gamma+\frac{1}{2})}U_{A,j}(s)\int_0^sU_{A,j}(u)^*
\bigl[\dot P_j(u),P_j(u)\bigr]U_T(u)\,dW_u.
\end{align*}
Taking the operator norm to both sides, then applying the triangle inequality pathwise and the expectation, we obtain
\begin{align}
\E\sup_{0\le s\le 1}\|U_{A,j}(s)-U_T(s)\|^2
&\le
\E\sup_{0\le s\le1}\Bigl[\frac{1}{T}\left\|X_j[\dot P_j](s)U_T(s)
-U_{A,j}(s)X_j[\dot P_j](0)\right\| \notag\\
&\quad
+\frac{1}{T}\left\|U_{A,j}(s)
\int_0^s U_{A,j}(u)^*
\Bigl([\dot P_j(u),P_j(u)]X_j[\dot P_j](u)
\right. \notag\\
&\quad\left.
-\frac{d}{du}X_j[\dot P_j](u)\Bigr)
U_T(u)\,du\right\|  \notag\\
&\quad
+\frac{\sigma^2\,T^{-2\gamma}}{2}\left\|U_{A,j}(s)
\int_0^sU_{A,j}(u)^*\Bigl[\bigl[\dot H_T(u),P_j(u)\bigr],P_j(u)\Bigr]
U_T(u)\,du\right\| \notag\\
&\quad
+\sigma\,T^{-(\gamma+\frac{1}{2})}\left\|U_{A,j}(s)
\int_0^sU_{A,j}(u)^*\bigl[\dot P_j(u),P_j(u)\bigr]
U_T(u)\,dW_u\right\|\Bigr]^2.
\label{eq:65}
\end{align}
By the triangle inequality and the submultiplicativity of the operator norm,
\[
\left\|X_j[\dot P_j](s)U_T(s)-U_{A,j}(s)X_j[\dot P_j](0) \right\| 
\le
\|X_j[\dot P_j](s)\|\,\|U_T(s)\|
+\|U_{A,j}(s)\|\,\|X_j[\dot P_j](0)\|.
\]
Therefore, by \eqref{eq:42},
\[
\left\|X_j[\dot P_j](s)U_T(s)-U_{A,j}(s)X_j[\dot P_j](0) \right\| 
\le
\|X_j[\dot P_j](s)\|
+\|X_j[\dot P_j](0)\|.
\]
From \eqref{eq:15},
\[
\left\|X_j[\dot P_j](s)U_T(s)-U_{A,j}(s)X_j[\dot P_j](0) \right\| 
\le
\frac{\|\dot P_j(s)\|}{g_j(s)}
+\frac{\|\dot P_j(0)\|}{g_j(0)}.
\]
Finally, using \eqref{eq:18},
\begin{equation}
\left\|X_j[\dot P_j](s)U_T(s)-U_{A,j}(s)X_j[\dot P_j](0) \right\| 
\le
\frac{\|\dot H_T(s)\|}{g_j(s)^2}
+\frac{\|\dot H_T(0)\|}{g_j(0)^2}
\qquad
\text{a.s..}
\label{eq:66}
\end{equation}
The triangle inequality and submultiplicativity yield
\begin{align}
&\left\|
U_{A,j}(s)\int_0^s
U_{A,j}(u)^*
\Bigl([\dot P_j(u),P_j(u)]X_j[\dot P_j](u)
-\frac{d}{du}X_j[\dot P_j](u)\Bigr)U_T(u)\,du\right\|\notag\\
\quad
&\le
\left\|U_{A,j}(s)\right\|\left\|\int_0^s
U_{A,j}(u)^*\Bigl([\dot P_j(u),P_j(u)]X_j[\dot P_j](u)
-\frac{d}{du}X_j[\dot P_j](u)\Bigr)U_T(u)\,du\right\|\notag\\
\quad
&\le
\left\|\int_0^sU_{A,j}(u)^*\Bigl([\dot P_j(u),P_j(u)]X_j[\dot P_j](u)
-\frac{d}{du}X_j[\dot P_j](u)\Bigr)U_T(u)\,du\right\|\notag\\
\quad
&\le
\int_0^s\left\|U_{A,j}(u)^*\Bigl([\dot P_j(u),P_j(u)]X_j[\dot P_j](u)
-\frac{d}{du}X_j[\dot P_j](u)\Bigr)U_T(u)\,\right\|du\notag\\
\quad
&\le
\int_0^s\left\| U_{A,j}(u)^* \right\| \left\|\Bigl([\dot P_j(u),P_j(u)]X_j[\dot P_j](u)
-\frac{d}{du}X_j[\dot P_j](u)\Bigr) \right\|\left\| U_T(u) \right\|\,du \notag\\
\quad
&\le
\int_0^s\left\|[\dot P_j(u),P_j(u)]\,X_j[\dot P_j](u)
-\frac{d}{du}X_j[\dot P_j](u)\right\|\,du\notag\\
\quad
&\le
\int_0^s\left(\left\|[\dot P_j(u),P_j(u)]\,X_j[\dot P_j](u)\right\|
+\left\|\frac{d}{du}X_j[\dot P_j](u)\right\|\right)\,du.
\label{eq:67}
\end{align}
Now, estimate the two terms in \eqref{eq:67} separately.

\medskip
\noindent
For the first term,
\[
\left\|
[\dot P_j(u),P_j(u)]\,X_j[\dot P_j](u)
\right\|
\le
2\,\|\dot P_j(u)\|\,\|X_j[\dot P_j](u)\|.
\]
By \eqref{eq:15},
\[
\left\|
[\dot P_j(u),P_j(u)]\,X_j[\dot P_j](u)
\right\|
\le
2\,\frac{\|\dot P_j(u)\|^2}{g_j(u)}.
\]
From \eqref{eq:18},
\begin{equation}
\left\|
[\dot P_j(u),P_j(u)]\,X_j[\dot P_j](u)
\right\|
\le
2\,\frac{\|\dot H_T(u)\|^2}{g_j(u)^3}.
\label{eq:68}
\end{equation}
For the second term, applying \eqref{eq:16} gives
\[
\left\|
\frac{d}{du}X_j[\dot P_j](u)
\right\|
\le
\frac{\|\ddot P_j(u)\|}{g_j(u)}
+
4\,\frac{\|\dot H_T(u)\|\,\|\dot P_j(u)\|}{g_j(u)^2}.
\]
Moreover, from \eqref{eq:18},
\[
\left\|
\frac{d}{du}X_j[\dot P_j](u)
\right\|
\le
\frac{\left\|
\frac{d}{du}[\dot P_j](u)
\right\|}{g_j(u)}
+4\,\frac{\|\dot H_T(u)\|^2}{g_j(u)^3}.
\]
Then, by \eqref{eq:10},
\begin{equation}
\left\|
\frac{d}{du}X_j[\dot P_j](u)
\right\|
\le
\frac{\left\|
\frac{d}{du}X_j[\dot H_T](u)
\right\|}{g_j(u)}
+4\,\frac{\|\dot H_T(u)\|^2}{g_j(u)^3}.
\label{eq:69}
\end{equation}
Substituting \eqref{eq:16} into \eqref{eq:69}, we obtain
\[
\left\|
\frac{d}{du}X_j[\dot P_j](u)
\right\|
\le
\frac{1}{g_j(u)}
\left(
\frac{\|\ddot H_T(u)\|}{g_j(u)}
+4\,\frac{\|\dot H_T(u)\|^2}{g_j(u)^2}
\right)
+4\,\frac{\|\dot H_T(u)\|^2}{g_j(u)^3}.
\]
Therefore,
\begin{equation}
\left\|
\frac{d}{du}X_j[\dot P_j](u)
\right\|
\le
\frac{\|\ddot H_T(u)\|}{g_j(u)^2}
+8\,\frac{\|\dot H_T(u)\|^2}{g_j(u)^3}.
\label{eq:70}
\end{equation}
Combine the bounds obtained in \eqref{eq:68} and \eqref{eq:70}, and substitute them into \eqref{eq:67} to obtain
\begin{align}
&\left\|
U_{A,j}(s)\int_0^s U_{A,j}^*(u)
\left([\dot P_j(u),P_j(u)]\,X_j[\dot P_j](u)-\frac{d}{du}X_j[\dot P_j](u)
\right)U_T(u)\,du\right\|\notag\\
\quad
&\le
\int_0^s
\left(\frac{\|\ddot H_T(u)\|}{g_j(u)^2}+10\,\frac{\|\dot H_T(u)\|^2}{g_j(u)^3}
\right)\,du
\qquad
\text{a.s..}
\label{eq:71}
\end{align}
By \eqref{eq:19}, \eqref{eq:20} and \eqref{eq:42}, this gives us,
\begin{align}
&\frac{\sigma^2\,T^{-2\gamma}}{2}\left\|U_{A,j}(s)
\int_0^sU_{A,j}(u)^*\Bigl[\bigl[\dot H_T(u),P_j(u)\bigr],P_j(u)\Bigr]
U_T(u)\,du\right\|\notag\\
\quad
&\le
\frac{\sigma^2\,T^{-2\gamma}}{2}\|U_{A,j}(s)\|\left\|\int_0^sU_{A,j}(u)^*
\Bigl[\bigl[\dot H_T(u),P_j(u)\bigr],P_j(u)\Bigr]U_T(u)\,du\right\|\notag\\
\quad
&\le
\frac{\sigma^2\,T^{-2\gamma}}{2}\int_0^s\|U_{A,j}(u)^*\|
\left\|\Bigl[\bigl[\dot H_T(u),P_j(u)\bigr],P_j(u)\Bigr]\right\|\|U_T(u)\|\,du\notag\\
\quad
&\le
\frac{\sigma^2\,T^{-2\gamma}}{2}\int_0^s\left\|\Bigl[\bigl[\dot H_T(u),P_j(u)\bigr],P_j(u)\Bigr]\right\|\,du\notag\\
\quad
&\le
\frac{\sigma^2\,T^{-2\gamma}}{2}\int_0^s\|\dot H_T(u)\|\,du
\qquad
\text{a.s..}
\label{eq:72}
\end{align}
From \eqref{eq:42}, we obtain
\begin{align}
&\sigma\,T^{-(\gamma+\frac{1}{2})}\left\|U_{A,j}(s)
\int_0^sU_{A,j}(u)^*\bigl[\dot P_j(u),P_j(u)\bigr]U_T(u)\,dW_u\right\|\notag\\
\quad
&\le
\sigma\,T^{-(\gamma+\frac{1}{2})}\|U_{A,j}(s)\|\left\|\int_0^sU_{A,j}(u)^*
\bigl[\dot P_j(u),P_j(u)\bigr]U_T(u)\,dW_u\right\|\notag\\
\quad
&\le
\sigma\,T^{-(\gamma+\frac{1}{2})}\left\|\int_0^sU_{A,j}(u)^*
\bigl[\dot P_j(u),P_j(u)\bigr]U_T(u)\,dW_u\right\|.
\label{eq:73}
\end{align}
Combine the bounds \eqref{eq:66}, \eqref{eq:71}, \eqref{eq:72} and \eqref{eq:73} into \eqref{eq:65}
\begin{align}
\E\sup_{0\le s\le 1}\|U_{A,j}(s)-U_T(s)\|^2
&\le
\E\sup_{0\le s\le 1}
\Bigl[\frac{C_j(s)}{T}
+\frac{\sigma^2\,T^{-2\gamma}}{2}\int_0^s\|\dot H_T(u)\|\,du
+\sigma\,T^{-(\gamma+\frac{1}{2})}\left\|M_s\right\|\Bigr]^2.
\label{eq:74}
\end{align}
The deterministic coefficient is
\[
C_j(s)
=
\frac{\|\dot H_T(s)\|}{g_j(s)^2}
+\frac{\|\dot H_T(0)\|}{g_j(0)^2}
+\int_0^s
\left(\frac{\|\ddot H_T(u)\|}{g_j(u)^2}+10\,\frac{\|\dot H_T(u)\|^2}{g_j(u)^3}
\right)\,du
\]
and
\[
M_s=\int_0^sU_{A,j}(u)^*\bigl[\dot P_j(u),P_j(u)\bigr]U_T(u)\,dW_u.
\]
Applying Young’s inequality to right side of \eqref{eq:74}, we obtain
\begin{align}
\E\sup_{0\le s\le 1}\|U_{A,j}(s)-U_T(s)\|^2
&\le
\E\sup_{0\le s\le 1}
\Biggl[\frac{3C_j(s)^2}{T^2}
+\frac{3\sigma^4\,T^{-4\gamma}}{4}
\Bigl(\int_0^s\|\dot H_T(u)\|\,du\Bigr)^2\notag\\
&\quad
+3\sigma^2\,T^{-(2\gamma+1)}\|M_s\|^2\Biggr]\notag\\
\quad
&\le
\frac{3}{T^2}\E\sup_{0\le s\le 1}C_j(s)^2
+\frac{3\sigma^4\,T^{-4\gamma}}{4}
\E\sup_{0\le s\le 1}\Bigl(\int_0^s\|\dot H_T(u)\|\,du\Bigr)^2\notag\\
&\quad
+3\sigma^2\,T^{-(2\gamma+1)}\E\sup_{0\le s\le 1}\|M_s\|^2\notag\\
\quad
&\le
\frac{3}{T^2}\sup_{0\le s\le 1}C_j(s)^2
+\frac{3\sigma^4\,T^{-4\gamma}}{4}
\sup_{0\le s\le 1}\Bigl(\int_0^s\|\dot H_T(u)\|\,du\Bigr)^2\notag\\
&\quad
+3\sigma^2\,T^{-(2\gamma+1)}\E\sup_{0\le s\le 1}\|M_s\|^2\notag\\
\quad
&\le
\frac{3}{T^2}\sup_{0\le s\le 1}C_j(s)^2
+\frac{3\sigma^4\,T^{-4\gamma}}{4}
\Bigl(\int_0^1\|\dot H_T(u)\|\,du\Bigr)^2\notag\\
&\quad
+3\sigma^2\,T^{-(2\gamma+1)}\E\sup_{0\le s\le 1}\|M_s\|^2.
\label{eq:75}
\end{align}
Let \(d=\dim\mathcal H\), and denote the corresponding Hilbert--Schmidt norm by
\(\|\cdot\|_{\mathrm{HS}}\). 

\medskip
\noindent
Take
\[
G_u=U_{A,j}(u)^*[\dot P_j(u),P_j(u)]U_T(u).
\]
Then,
\[
M_s=\int_0^s G_u\,dW_u.
\]

\medskip
\noindent
By \eqref{eq:42} and \eqref{eq:17},
\begin{align}
\|G_u\|_{\mathrm{HS}}^2
&\le d\,\|G_u\|^2\notag\\
\quad
&= d\,\|[\dot P_j(u),P_j(u)]\|^2\notag\\
\quad
&= d\,\|\dot P_j(u)\|^2
\qquad\text{a.s..}
\label{eq:76}
\end{align}
Since \(\dot P_j\) is deterministic and continuous on \([0,1]\),
 \[
 \E\int_0^1\|G_u\|_{\mathrm{HS}}^2\,du<\infty.
 \]
Thus, \(\{M_s\}_{s \in [0,1]}\) is a square-integrable martingale.

\medskip
\noindent
Applying Doob's \(L^2\)-maximal inequality to the submartingale
\(\|M_s\|_{\mathrm{HS}}\) then using Itô's isometry~\cite[Theorem~4.3.5]{Kuo2006} entrywise and \eqref{eq:76}, we obtain
\begin{align*}
\E\!\left[\sup_{0\le s\le1}\|M_s\|^2\right]
&\le \E\!\left[\sup_{0\le s\le1}\|M_s\|_{\mathrm{HS}}^2\right] \\
&\le 4\,\E\|M_1\|_{\mathrm{HS}}^2 \\
&= 4\,\E\int_0^1\|G_u\|_{\mathrm{HS}}^2\,du \\
&\le 4d\int_0^1\|\dot P_j(u)\|^2\,du.
\end{align*}
From \eqref{eq:18},
\begin{align}
\E\sup_{0\le s\le 1}
\left\|\int_0^sU_{A,j}(u)^*
\bigl[\dot P_j(u),P_j(u)\bigr]U_T(u)\,dW_u\right\|^2
&\le
4d\int_0^1\frac{\|\dot H_T(u)\|^2}{g_j(u)^2}\ du.
\label{eq:77}
\end{align}
By \eqref{eq:75} and \eqref{eq:77},
\begin{align}
\E\sup_{0\le s\le 1}\|U_{A,j}(s)-U_T(s)\|^2
&\le
\frac{3}{T^2}\sup_{0\le s\le 1}C_j(s)^2
+\frac{3\sigma^4\,T^{-4\gamma}}{4}
\Bigl(\int_0^1\|\dot H_T(u)\|\,du\Bigr)^2\notag\\
&\quad
+12\sigma^2\,d\,T^{-(2\gamma+1)}
\Bigl(\int_0^1\frac{\|\dot H_T(u)\|^2}{g_j(u)^2}\ du\Bigr).
\label{eq:78}
\end{align}
Since $P_j(0)$ is an orthogonal projection,
\[
\|P_j(0)\|=1.
\]
Therefore,
\begin{align}
\|(U_{A,j}(s)-U_T(s))P_j(0)\|
&\le
\|U_{A,j}(s)-U_T(s)\|\,\|P_j(0)\| \notag\\
&\le
\|U_{A,j}(s)-U_T(s)\|.
\label{eq:79}
\end{align}
From \eqref{eq:78} and \eqref{eq:79},
\begin{align}
\E\sup_{0\le s\le 1}\|(U_{A,j}(s)-U_T(s))P_j(0)\|^2
&\le
\frac{3}{T^2}\sup_{0\le s\le 1}C_j(s)^2
+\frac{3\sigma^4\,T^{-4\gamma}}{4}
\Bigl(\int_0^1\|\dot H_T(u)\|\,du\Bigr)^2\notag\\
&\quad
+12\sigma^2\,d\,T^{-(2\gamma+1)}
\int_0^1\frac{\|\dot H_T(u)\|^2}{g_j(u)^2}\ du.
\label{eq:80}
\end{align}
\end{proof}
\noindent
The final estimate bounds the maximal mean-square error
\[
\mathcal E_T
=
\E\sup_{0\le s\le 1}\|(U_{A,j}(s)-U_T(s))P_j(0)\|^2.
\]
There exist finite nonnegative constants \(K_0,K_1,K_2\), independent of \(T\), such that
\begin{align*}
\sup_{0\le s\le1}C_j(s)^2
&\le 
K_0,\notag\\
\quad
\frac{1}{4}\Bigl(\int_0^1\|\dot H_T(u)\|\,du\Bigr)^2
&\le 
K_1,\notag\\
\quad
d\int_0^1\frac{\|\dot H_T(u)\|^2}{g_j(u)^2}\,du
&\le 
K_2.
\end{align*}
Then \eqref{eq:80} becomes
\begin{equation*}
\mathcal E_T
\le
3K_0T^{-2}+3K_1\sigma^4T^{-4\gamma}+12K_2\sigma^2T^{-(2\gamma+1)}.
\end{equation*}
Consequently,
\begin{equation*}
\mathcal E_T
\le
C_0T^{-2}+C_1\sigma^4T^{-4\gamma}+C_2\sigma^2T^{-(2\gamma+1)},
\end{equation*}
with
\[
C_0=3K_0,\qquad
C_1=3K_1,\qquad
C_2=12K_2.
\]
The first, second, and third terms arise from the deterministic coefficient, the Itô correction, and the martingale, respectively.

\medskip
\noindent
For fixed \(\sigma>0\), the estimate yields the following asymptotic upper bounds as 
\(T\to\infty\):
\[
\mathcal E_T=
\begin{cases}
O(1),
& \gamma=0:
\quad\text{fixed nonzero physical noise amplitude; convergence remains unresolved},
\\[6pt]
O(T^{-4\gamma}),
& 0<\gamma<\tfrac12:
\quad\text{the second term decays most slowly in the upper bound},
\\[6pt]
O(T^{-2}),
& \gamma\geq\tfrac12:
\quad\text{the bound retains the deterministic order}.
\end{cases}
\]
Equivalently, the root-mean-square maximal error satisfies
\[
\mathcal E_T^{1/2}=
\begin{cases}
O(1), & \gamma=0,\\[4pt]
O(T^{-2\gamma}), & 0<\gamma<\tfrac12,\\[4pt]
O(T^{-1}), & \gamma\geq\tfrac12.
\end{cases}
\]
Every $\gamma>0$ therefore guarantees convergence, while \(\gamma\ge\frac{1}{2}\) recovers the deterministic \(O(T^{-1})\) upper bound for the root-mean-square maximal error. For \(\gamma=0\), convergence is not established by this estimate. The \(O(1)\) upper bound does not establish a nonzero limiting error.

\section{Discussion}
The main theorem bounds the maximal mean-square error between the exact stochastic propagator \(U_T(s)\) and stochastic adiabatic intertwiner \(U_{A,j}(s)\) on the initial eigenspace corresponding to the selected isolated eigenvalue. The same real-valued Wiener process \(W_s\) drives both operators, so the adiabatic evolution includes the random dynamical phase. Taking the supremum before the expectation captures temporary deviations that may be missed by considering only the endpoint.

The spectral gap remains central. The deterministic coefficient \(C_j(s)\) in \eqref{eq:80} involves first and second derivatives of the Hamiltonian \(H_T(s)\) and inverse powers of the spectral gap \(g_j(s)\), while the third term in \eqref{eq:80} depends on the integral of \(\frac{\|\dot H_T\|^2}{g_j^2}\). Therefore, small gaps weaken both estimates. For the second term in \eqref{eq:80}, the double commutator identity bounds its coefficient through the integrated variation of the Hamiltonian without an explicit inverse gap factor.

For a fixed Hamiltonian path, take \(\alpha=\sigma T^{-\gamma}\) with constant \(\sigma>0\) and \(\gamma>0\). Every \(\gamma>0\) makes all terms vanish as \(T\to\infty\), establishing maximal mean-square convergence, provided the path-dependent coefficients remain bounded independently of \(T\). For \(0<\gamma<\tfrac{1}{2}\), the second term decays most slowly, giving an \(O(T^{-4\gamma})\) bound. At \(\gamma=\frac{1}{2}\), all three mean-square terms scale as \(T^{-2}\), giving an \(O(T^{-2})\) maximal mean-square error \(\mathcal E_T\) and an \(O(T^{-1})\) root-mean-square maximal error \(\mathcal E_T^{1/2}\). This recovers the deterministic rate while the scalar noise coefficient in rescaled time \(s\) remains \(\sigma\). Thus, for \(\gamma>\tfrac12\), the second and third contribution decay faster than the first term in the maximal mean-square error \(\mathcal E_T\).

The scalar noise coefficient in rescaled time \(s\) is \(\sigma\sqrt{T}\) for fixed nonzero physical noise amplitude \(\gamma=0\). The first term decays as \(T^{-2}\), while the second term contributes a term proportional to \(\sigma^4\) with no inverse power of the total evolution time \(T\) in the maximal mean-square error \(\mathcal E_T\). This term gives no positive lower bound on the error. Therefore, the estimate establishes neither maximal mean-square convergence nor a breakdown of adiabaticity.

The result provides convergence for every \(\gamma>0\) under the stated uniform bounds and noise scaling \(\alpha\) with fixed \(\sigma\). We bound the terms in the difference between the exact stochastic propagator \(U_T(s)\) and stochastic adiabatic intertwiner \(U_{A,j}(s)\) separately. Therefore, we may miss cancellations that would give a faster error decay. For fixed nonzero noise physical noise amplitude \((\alpha=\sigma>0)\), the upper bound remains inconclusive as to whether the error tends to zero as \(T\to\infty\).

To extend the intertwining proof to several independent noise sources, we require each coupling operator to commute with the spectral projection \(P_j(s)\) for every \(s\in[0,1]\). A noncommuting coupling can move a state out of the selected eigenspace. The off-diagonal blocks of the averaged density matrix could show whether coherence between eigenspaces decays during evolution. If the spectral gap closes, then the inverse powers of the spectral gap in our bound can diverge. In infinite dimensions, the unitarity proof must show that each propagator is norm-preserving and that its range is the entire Hilbert space.

\section*{Disclosure statement}
The authors declare that they have no competing interests.

\section*{Declaration of AI utilization}
We utilized ChatGPT (OpenAI; GPT-6 Sol and GPT-6 Astra) and Gemini (Google; Gemini 3 Flash and Gemini 3 Pro) to explore research ideas, edit the manuscript for clarity, and organize references. For the fourth paragraph of the Introduction, both tools helped us identify publications on stochastic Schrödinger equations and the quantum adiabatic theorem. We selected the publications to cite and checked our descriptions against the original sources. We reviewed terms of use for ChatGPT (OpenAI, GPT-6 Sol and GPT-6 Astra) and Gemini (Google, Gemini 3 Flash and Gemini 3 Pro) and confirmed that our use was suitable for publication. We independently verified the mathematical arguments and references. We confirm the originality and accuracy of the final manuscript and take full responsibility for its content.

\end{document}